\documentclass[%
	aps,
	prd,
	reprint,
	a4paper,
	floatfix,
	nofootinbib,
	superscriptaddress
]
{revtex4-2}

\usepackage[utf8]{inputenc}
\usepackage[T1]{fontenc}
\usepackage[american]{babel}
\usepackage{graphicx}
\usepackage{xcolor}
\usepackage{mathtools}
\usepackage{amssymb}
\usepackage{lmodern}
\usepackage{booktabs}
\usepackage{enumitem}
\usepackage{orcidlink}
\usepackage[capitalize]{cleveref}
\usepackage{microtype}

\newcommand{\mytitle}{Curve fitting on a quantum annealer for an advanced navigation method}

\definecolor{refcolor}{cmyk}{0.9 0.9 0 0}

\hypersetup{%
	pdftitle = {\mytitle},
	pdfauthor = {Philipp Isserstedt, Daniel Jaroszewski, Wolfgang Mergenthaler, Felix Paul, Bastian Harrach},
	bookmarksopen = true,
	colorlinks = true,
	allcolors = refcolor
}

\newcommand*{\FCE}{%
	FCE Frankfurt Consulting Engineers GmbH, %
	Bessie-Coleman-Str.~7, %
	60549 Frankfurt am Main, %
	Germany%
}

\newcommand*{\AN}{%
	Anaqor AG, %
	Keithstr.~6, %
	10787 Berlin, %
	Germany%
}

\newcommand*{\GU}{%
	Institut f\"{u}r Mathematik, %
	Goethe-Universit\"{a}t Frankfurt, %
	Robert-Mayer-Str.~10, %
	60325 Frankfurt am Main, %
	Germany%
}

\newcommand*{\DLR}{%
	DLR-Institut f\"{u}r Weltraumforschung, %
	Rutherfordstr.~2, %
	12489 Berlin, %
	Germany%
}

\newcommand*{\+}{\hspace*{0.08335em}}
\newcommand*{\?}{\hspace*{-0.08335em}}
\newcommand*{\opt}{\ast}
\newcommand*{\calF}{\mathcal{F}}
\newcommand*{\calN}{\mathcal{N}}
\newcommand*{\setB}{\{0,1\}}
\newcommand*{\setN}{\mathbb{N}}
\newcommand*{\setR}{\mathbb{R}}

\newcommand*{\setRppos}{\setR_{>0}}
\newcommand*{\mystate}[1]{x_{#1}}
\newcommand*{\action}[1]{u_{#1}}
\newcommand*{\setStates}[1]{X_{#1}}
\newcommand*{\setActions}[1]{U_{#1}}
\newcommand*{\valFuncSymbol}{V}
\newcommand*{\valFunc}[1]{\valFuncSymbol_{#1}}
\newcommand*{\valFuncOpt}[1]{\valFunc{#1}^{\opt}}
\newcommand*{\valFuncOptAppr}[1]{\hat{\valFuncSymbol}_{#1}^{\opt}}
\newcommand*{\transFunc}[1]{\Gamma_{#1}}
\newcommand*{\costFunc}[1]{C_{#1}}
\newcommand*{\vmax}{v_{\textup{max}}}
\newcommand*{\argmin}{\operatorname*{arg\,min}}
\newcommand*{\grad}{\nabla}

\newcommand*{\diag}{\operatorname{diag}}
\newcommand*{\triang}{\Lambda}
\newcommand*{\ZIsing}{Z_{\textup{Ising}}}
\newcommand*{\ZQUBO}{Z_{\textup{QUBO}}}

\newcommand*{\ape}{\operatorname{ape}}

\newcommand*{\rmse}{\operatorname{rmse}}
\newcommand*{\foptcl}{f_{\textup{cl}}^{\opt}}
\newcommand*{\fopttb}{f_{\textup{tb}}^{\opt}}
\newcommand*{\foptqa}{f_{\textup{qa}}^{\opt}}
\newcommand*{\supp}{\operatorname{supp}}
\crefname{table}{Tab.}{Tabs.}
\Crefname{table}{Table}{Tables}

\AtBeginDocument{\abovedisplayshortskip=0em}

\begin{document}

\title{\mytitle}

\author{Philipp Isserstedt\,\orcidlink{0000-0003-1989-0886}}
\email{philipp.isserstedt@frankfurtconsultingengineers.de}
\affiliation{\FCE}

\author{Daniel Jaroszewski}
\email{daniel.jaroszewski@frankfurtconsultingengineers.de}
\affiliation{\FCE}
\affiliation{\GU}

\author{Wolfgang Mergenthaler}
\email{wolfgang.mergenthaler@frankfurtconsultingengineers.de}
\affiliation{\FCE}

\author{Felix Paul}
\email{felix.paul@dlr.de}
\affiliation{\AN}
\affiliation{\DLR}

\author{Bastian Harrach\,\orcidlink{0000-0002-8666-7791}\,}
\email{harrach@math.uni-frankfurt.de}
\affiliation{\GU}


\begin{abstract}
We explore the applicability of quantum annealing to the approximation task of curve fitting.
To this end, we consider a function that shall approximate a given set of data points and is
written as a finite linear combination of standardized functions, e.g., orthogonal polynomials.
Consequently, the decision variables subject to optimization are the coefficients of that expansion.
Although this task can be accomplished classically, it can also be formulated as a quadratic
unconstrained binary optimization problem, which is suited to be solved with quantum annealing.
Given the size of the problem stays below a certain threshold, we find that quantum annealing yields
comparable results to the classical solution. Regarding a real-world use case, we discuss the
problem to find an optimized speed profile for a vessel using the framework of dynamic programming
and outline how the aforementioned approximation task can be put into play. Similar to the curve
fitting task, our findings indicate that quantum annealing is currently only feasible if the
routing problem is modeled sufficiently small and sparse.
\\[0.75em]
\end{abstract}

\maketitle

\section{Introduction}
\label{introduction}

In 2020, the International Maritime Organization issued a recommendation to
navigate just in time such that fuel consumption is minimized while arrival
at the destination is guaranteed to be on time \cite{IMO:2020a}. In this case,
the approximation corresponds to an optimal speed profile for a vessel along
a given trajectory in the sense that fuel consumption is minimized while
arrival at the destination is guaranteed to be just in time.
The given task is in fact an element of variational calculus and functional
analysis. If one discretizes the problem in an attempt to generate a
discrete sequence of speed values to support those objectives, one quickly ends
up in the domain of dynamic programming and finds the Bellman equation to be
a highly valuable tool \cite{Bellman:1952a,Bellman:1953a,Bellman:1954a,Bellman:1957a,%
Neumann:2002a,Bertsekas:2017a,Bertsekas:2012a,Sutton:2018a,Bertsekas:2019a}.

Using dynamic programming to address large real-world use cases numerically
necessitates a discretization of the problem's state space. As a consequence, this
entails an error in the accuracy of the so-called value function, which is the
solution of the Bellman equation, and it is common practice to employ approximation
techniques; see, e.g., Refs.~\cite{Bertsekas:2019a,Lagoudakis:2010} and references
therein.

One of the most fundamental approximation tasks is curve fitting, where a function
is sought that approximates a given set of data, i.e., the difference between the
observed data and the approximation shall become minimal. In this work, we formulate
curve fitting as a quadratic unconstrained binary optimization (QUBO)
problem \cite{Kochenberger:2014a,Glover:2018a} that is suited to be
solved with quantum computing---particularly quantum annealing
\cite{Apolloni:1989a,Kadowaki:1998a,Brooke:1999a,Farhi:2000a,Farhi:2001a,%
Farhi:2002a,Aharonov:2004a,Somma:2012a,Lucas:2013a,Bian:2014a,McGeoch:2014a,%
Albash:2016a,Venegas-Andraca:2018a,Yarkoni:2021a}.

Furthermore, in this context we discuss the task of finding an optimized speed
profile for a vessel in the sense of minimal fuel consumption while arriving just
in time using the framework of dynamic programming. Having the voyage of a vessel
in mind, it might very well be the case that finding an optimized speed profile
requires a lot of support points in order to densely discretize the state space,
especially as soon as the latter becomes high-dimensional, and approximations are
practically inevitable. Here, we use the QUBO formulation of curve fitting solved
on a quantum annealer to approximate the value function with respect to its
dependency on the state space variable.

With quantum computing currently being an emerging technology still in its
infancy, we do not expect to outperform a conventional computer. Particularly in
view of the fact that the only current commercially available quantum annealing
devices manufactured by D-Wave \cite{Berkeley:2009a,Johnson:2009a,Harris:2009a,%
Harris:2010a,Johnson:2011a,Dickson:2013a,Lanting:2013a,Lanting:2014a,Bunyk:2014a,%
McGeoch:2014a,McGeoch:2019a} are not full-fledged enough yet, and it is a priori
not clear whether their use will result in practical advantages for a specific
problem \cite{Parekh:2016a,Venegas-Andraca:2018a,Juenger:2021a,Yarkoni:2021a}.
Other application- and industry-focused developement toward quantum annealing
devices is carried out by, e.g., the NEC Corporation \cite{Lin:2014,NEC:2022,NEC:2023}
and Qilimanjaro Quantum Tech \cite{Canivell:2021,Palacios:2024}. There are also
quantum annealing-inspired devices tailored to solve QUBO problems, e.g., the coherent
Ising machine \cite{Inagaki:2016,McMahon:2016}.

Nevertheless, in our opinion it is worth to explore to what extent quantum annealing
can be applied to real-world use cases that arise in an application- and industry-focused
environment (see Ref.~\cite{Yarkoni:2021a} for an overview) and to judge the quality of
the solutions compared to classical methods. This work lies within that scope.

The remainder of this article is organized as follows. In \cref{curve_fitting},
we formulate curve fitting as a QUBO and present our numerical results in
\cref{numerical_experiments}. After that, \cref{speed_profile} focuses on
finding an optimized speed profile using dynamic programming, and we finally
conclude in \cref{summary}.

\section{Curve fitting as a QUBO}
\label{curve_fitting}

Let $X, Y \subset \setR$, $n \in \setN$, and
$\{ (x_{i}, y_{i}) \}_{i=0}^{n-1} \subset X \times Y$ the observed data points.
The task is to find an optimal fit function $f^{\opt} \colon X \to Y$ in the sense that
the sum of squares of the differences between the observed data $y_{i}$ and approximated
values $f^{\opt}(x_{i})$ is minimized (method of least squares) \cite{Quarteroni:2007a}.
With $\calF$ denoting a suitable function space, e.g.,
$\calF = \{ f \colon \setR \to \setR \mid f \text{~continuous} \}$,
we are concerned with the optimization problem
\begin{equation}
	\label{eq:least_squares}
	\min_{f \in \calF} \sum_{i=0}^{n-1} (y_{i} - f(x_{i}))^{2} \, ,
\end{equation}
whose solution is the desired optimal fit $f^\opt$.%
\footnote{By construction of minimizing residuals, the least squares method always
yields a solution, i.e., an optimal fit $f^{\opt}$ always exists. In our case, where
the fit function is written as a finite linear combination of standardized functions according
to \cref{eq:function_expansion}, the solution is unique if the matrix $W$, defined
later in \cref{eq:classical_01b}, has full rank.}
For example, if $f(x) = a_{0} + a_{1} x$, the minimization task \labelcref{eq:least_squares}
corresponds to identifying the optimal parameters $a_{0}^{\opt}, a_{1}^{\opt} \in \setR$
for an ordinary simple linear regression model.

Here, we allow for a more general expression for $f$ and express it as a finite
linear combination of standardized functions
$\{ \phi_{j} \colon X \to Y \}_{j=0}^{m-1} \subset \calF$, $m \in \setN$, i.e.,
\begin{equation}
	\label{eq:function_expansion}
	f(x) = \sum_{j=0}^{m-1} c_{j} \phi_{j}(x)
\end{equation}
with $c_{j} \in \setR$. Finding an optimal fit $f^{\opt}$ thus translates into finding
optimal coefficients $c_{j}^{\opt}$ (for given $\phi_{j}$). For example, such expansions
are frequently employed within a certain area of high-energy physics where they aid to handle
the numerical complexity of the occurring equations \cite{Eichmann:2016yit,Sanchis-Alepuz:2017jjd}.
Two customary choices for $\phi_{j}$ are:
\begin{itemize}[leftmargin=*, widest=(ii)]
	\item[(i)] orthogonal polynomials such as Chebyshev polynomials. The first
	kind of the latter, $T_{j}$, are defined via the recurrence relation
	\begin{equation}
		\label{eq:chebyshev}
		T_{j}(x) = 2x T_{j-1}(x) - T_{j-2}(x)
	\end{equation}
	for $j \in \setN_{>1}$ with $T_{1}(x) = x$ and $T_{0}(x) = 1$ \cite{Olver:2010a};
	\item[(ii)] triangular functions $\triang_{j}$ that yield a piecewise-linear
	approximation. They are defined with respect to supporting points
	$\tilde{x}_{0} < \dotso < \tilde{x}_{m-1}$, $m \geq 2$, with $\tilde{x}_{0} = x_{0}$
	and $\tilde{x}_{m-1} = x_{n-1}$ according to
	\begin{equation}
		\label{eq:triangular}
		\triang_{j}(x)
		=
		\begin{dcases}
			\frac{x - \tilde{x}_{j-1}}{\tilde{x}_{j} - \tilde{x}_{j-1}}
			& \text{if~} x \in [\tilde{x}_{j-1}, \tilde{x}_{j})
			\\[0.25em]
			\frac{\tilde{x}_{j+1} - x}{\tilde{x}_{j+1} - \tilde{x}_{j}}
			& \text{if~} x \in [\tilde{x}_{j}, \tilde{x}_{j+1})
			\\[0.25em]
			0 & \text{otherwise}
		\end{dcases}
	\end{equation}
	for $j = 0, \dotsc, m - 1$, where it is understood that only the first (second)
	branch defines the nonzero values of the last (first) triangular function
	$\triang_{j=m-1}$ ($\triang_{j=0}$).
\end{itemize}

\subsection{Classical solution}
\label{curve_fitting_classical}

We aim to determine the expansion coefficients $c_{j}$ in
\cref{eq:function_expansion} according to the least-squares
principle as given by the optimization problem \labelcref{eq:least_squares}.
With the abbreviations (where $j, k = 0, \dotsc, m-1$)
\begin{align}
	\label{eq:classical_01a}
	b_{j} &= \sum_{i=0}^{n-1} y_{i} \phi_{j}(x_{i}) \, ,
	\\[0.25em]
	\label{eq:classical_01b}
	W_{jk} &= \sum_{i=0}^{n-1} \phi_{j}(x_{i}) \phi_k(x_{i}) \, ,
\end{align}
vectors $c = [c_{0}, \dotsc, c_{m-1}]^{\top} \in \setR^{m}$,
$b = [b_{0}, \dotsc, b_{m-1}]^{\top} \in \setR^{m}$,
and the symmetric matrix $W \in \setR^{m \times m}$ with entries
$W_{jk}$, the optimization problem upon inserting \cref{eq:function_expansion}
into \labelcref{eq:least_squares} now takes the form
\begin{equation}
	\label{eq:classical_02}
	\begin{dcases}
		\min_{c \in \setR^{m}} Z(c) \, ,
		\\
		Z(c) = c^{\top} W c - 2 c^{\top} b \, .
	\end{dcases}
\end{equation}
The optimal weight coefficient vector $c^{\opt}$ thus follows from
\begin{equation}
	\label{eq:classical_03}
	\grad Z(c^{\opt}) = 2 \, (W c^{\opt} - b) \stackrel{!}{=} 0
\end{equation}
and is explicitly given by
\begin{equation}
	\label{eq:classical_04}
	c^{\opt} = W^{-1} b
\end{equation}
provided $W$ is invertible. If this is not the case, in practice one could
resort to its Moore--Penrose pseudoinverse \cite{Ben-Israel:2003a}.

\subsection{Quantum annealing-oriented approach}
\label{curve_fitting_qubo}

Starting from the objective function as given in \labelcref{eq:classical_02}
in its expanded form
\begin{equation}
	\label{eq:quantum_01}
	Z(c)
	=
	\sum_{j=0}^{m-1} \sum_{k=0}^{m-1} c_{j} W_{jk} c_{k}
	-
	2 \sum_{j=0}^{m-1} c_{j} b_{j} \, ,
\end{equation}
we have a quadratic function in the components of $c$. In order to express the
optimization problem \labelcref{eq:classical_02} as a QUBO, whose decision
variables' domain is the set $\setB$, we write each coefficient $c_{j}$ in binary
fixed-point format and use two's complement to represent negative values.
This yields the following expression for the coefficients:
\begin{equation}
	\label{eq:quantum_02a}
	c_{j} = \sum_{r=0}^{d-1} \psi_{jr} \sigma_{r} 2^{r-p} \, ,
\end{equation}
where $\psi_{jr} \in \setB$, $d \in \setN$ defines how many binary digits are used,
$p = 0, \dotsc, d - 1$ indicates the location of the fixed decimal point, and
\begin{equation}
	\label{eq:quantum_02b}
	\sigma_{r} =
	\begin{dcases}
		-1 & \text{if~} r = d - 1
		\\[0.25em]
		1 & \text{otherwise}
	\end{dcases}
\end{equation}
accounts for two's complement, i.e., the most significant bit enters the sum with a
negative sign. For example, for $p = 0$, \cref{eq:quantum_02a} allows to represent
all integers between $-(2^{d-1})$ and $2^{d-1} - 1$. The quantities $\psi_{jr}$ are then
the new decision variables represented by the qubit states of a quantum annealing device.

In principle, one could choose different values for $d$ and $p$ for each coefficient
$c_{j}$ matching the individual order of magnitude and location of the decimal point of each
coefficient. In practice, however, these are (usually) not known a priori. Thus, and for the
sake of simplicity, we choose the same $d$ and $p$ for all coefficients.

Now, inserting \cref{eq:quantum_02a} into \labelcref{eq:quantum_01} yields a new
objective function depending on $\psi = [\psi_{00}, \psi_{01}, \dotsc]^\top \in \setB^{md}$,
i.e., on all coefficients $\psi_{jr}$ of the binary representation of the expansion coefficients; to wit
\begin{align}
	\label{eq:quantum_03}
	Z(\psi)
	&=
	\sum_{j,k=0}^{m-1} \sum_{r,s=0}^{d-1} \psi_{jr} \psi_{ks} \sigma_{r} \sigma_{s} W_{jk} 2^{r+s-2p}
	\notag\\
	&\hphantom{=\;}-
	\sum_{j=0}^{m-1} \sum_{r=0}^{d-1} \psi_{jr} \sigma_{r} b_{j} 2^{r-p+1} \, .
\end{align}
This equation can also we written as
\begin{equation}
	\label{eq:quantum_04}
	Z(\psi)
	=
	\sum_{\mu,\nu=0}^{md-1} \psi_{\mu} W_{\mu\nu}' \psi_{\nu}
	-
	\sum_{\mu=0}^{md-1} \psi_{\mu} b_{\mu}' \, ,
\end{equation}
where we introduced $\psi_{\mu} = \psi_{\hat{\mu}\bar{\mu}}$,
$b_{\mu}' = 2^{\bar{\mu}-p+1} \sigma_{\bar{\mu}} b_{\hat{\mu}}$, and
$W_{\mu\nu}' = 2^{\bar{\mu}+\bar{\nu}-2p} \sigma_{\bar{\mu}} \sigma_{\bar{\nu}} W_{\hat{\mu}\hat{\nu}}$
with $\hat{\mu} = \lfloor \mu / d \rfloor$, $\bar{\mu} = \mu \bmod d$ (analogously for $\nu$).
$b_{\mu}'$ and $W_{\mu\nu}'$ are the entries of the quantities $b' \in \setR^{md}$ and
$W' \in \setR^{(md) \times (md)}$, respectively. These definitions allow us to write the objective
function in a rather compact form:
\begin{equation}
	\label{eq:quantum_05}
	Z(\psi) = \psi^{\top} W' \psi - \psi^{\top} b' \, .
\end{equation}

Finally, in order to cast this objective function into the canonical QUBO form,
we exploit that the entries of the binary vector $\psi$ are idempotent,
$\psi_{\mu}^{2} = \psi_{\mu}$ for all $\mu = 0, \dotsc, md - 1$. The linear term
$\psi^{\top} b'$ can therefore be written as a contribution that is subtracted from
the diagonal elements of $W'$. Thus, our final expression for the least-squares optimization
problem \labelcref{eq:least_squares} written as a QUBO reads
\begin{equation}
	\label{eq:quantum_06}
	\min_{\psi \in \setB^{md}} \psi^{\top} Q \psi
\end{equation}
with the QUBO matrix
\begin{equation}
	\label{eq:quantum_07}
	Q = W' - \diag(b_{0}', \dotsc, b_{md-1}') \, .
\end{equation}
In a last step, the solution $\psi^{\opt}$ of the QUBO problem is plugged into
\cref{eq:quantum_02a} to obtain the actual, real-valued optimal coefficients
$c_{j}^{\opt} = \sum_{r} \psi_{jr}^{\opt} \sigma_{r} 2^{r-p}$.

Apparently, the number of required qubits, i.e., the size of the QUBO problem, depends only
on the number of standardized functions $m$ used for the expansion of the fit function, see
\cref{eq:function_expansion}, and the number of digits $d$ used for the binary
representation of the expansion coefficients:
\begin{equation}
	Q \in \setR^{(md) \times (md)} \, .
\end{equation}
In particular, the size of the QUBO matrix is independent of the number of the
to-be-approximated data points $n$ (see \cref{eq:least_squares}); the latter enters only
implicitly in the course of calculating the entries of $Q$
(see Eqs.~\labelcref{eq:classical_01a} and \labelcref{eq:classical_01b}).

\section{Numerical experiments}
\label{numerical_experiments}

In the following, we present and discuss our numerical results of the curve fitting
task by means of a QUBO problem as described above.

Since we shall compare results obtained on a classical computer using a tabu search
\cite{Glover:1989a,Glover:1990a} with results from a D-Wave quantum annealing
device, we sketch very briefly the idea behind the latter (see, e.g.,
Refs.~\cite{Bunyk:2014a,McGeoch:2014a,McGeoch:2019a,DWave-Doc} and references therein
for more details). Based on the concept of adiabatic quantum computing \cite{Albash:2016a}
and relaxing some of its properties/conditions to facilitate a feasible
realization, quantum annealing is a metaheuristic method \cite{McGeoch:2014a} to find
low-energy solutions---including the ground state(s)---of an Ising spin model
\cite{Ising:1925a,Itzykson:1989a}. For a system of spins
$s = [s_{0}, \dotsc, s_{N-1}]^{\top} \in \{ -1, 1 \}^{N}$, $N \in \setN$, a ground
state is a spin configuration that globally minimizes the objective function
\begin{equation}
	\label{eq:qann_01}
	\ZIsing(s)
	=
	\sum_{i=0}^{N-1} \sum_{j=i+1}^{N-1} J_{ij} s_{i} s_{j}
	+
	\sum_{i=0}^{N-1} h_{i} s_{i} \, ,
\end{equation}
where $J_{ij} \in \setR$ and $h_{i} \in \setR$ denote the nearest-neighbor interaction
between the $i\+$\textsuperscript{th} and $j\+$\textsuperscript{th} spin (coupling strength)
and an external magnetic field acting on the $i\+$\textsuperscript{th} spin (qubit bias),
respectively.

The connection to the corresponding QUBO model with objective function
$\ZQUBO(z) = z^{\top} R \+ z$ with $z \in \setB^{N}$ and an upper-triangular matrix
$R \in \setR^{N \times N}$ is established by the transformation $s_{i} = 2 z_{i} - 1$
for $i = 0, \dotsc, N - 1$. The off- and on-diagonal entries of $R$ are appropriately
mapped onto the coupling strengths and qubit biases, respectively. It is then straightforward
to show that the Ising and QUBO objective functions differ only by an additive constant
in the sense that
\begin{equation}
	\label{eq:qann_02}
	\ZIsing(s) \bigr\rvert_{s_{i}=2z_{i}-1} = \ZQUBO(z) + \textup{const.} \, ,
\end{equation}
i.e., an approximate or exact solution of a given QUBO problem can be equivalently
obtained by finding low-energy or ground states of an Ising model, respectively.

\subsection{Sample data}
\label{sample_data}

Before we continue with our actual results, we define the sample data
$\{ (x_{i}, y_{i}) \}_{i=0}^{n-1}$ (see \cref{curve_fitting}) that shall be used for
our numerical experiments. The abscissa values are evenly spaced over the unit interval
according to
\begin{equation}
	\label{eq:sample_data_01}
	x_{i} = \frac{i}{n-1} \in [0,1] \, .
\end{equation}
For the ordinates, we choose the following four different types:
\begin{align}
	\label{eq:sample_data_02a}
	\textup{linear:}
	\quad
	y_{i} &= \frac{1}{2} x_{i} + 1 + \varepsilon_{i} \, ,
	\\[0.5em]
	\label{eq:sample_data_02b}
	\textup{quadratic:}
	\quad
	y_{i} &= \frac{3}{4} x_{i}^{2} + \varepsilon_{i} \, ,
	\\[0.5em]
	\label{eq:sample_data_02c}
	\textup{cubic:}
	\quad
	y_{i} &= \frac{3}{4} x_{i}^{3} + \frac{1}{4} x_{i} + \varepsilon_{i} \, ,
	\\[0.5em]
	\label{eq:sample_data_02d}
	\textup{trigonometric:}
	\quad
	y_{i} &= \sin(2 \pi x_{i}) \cos(2 \pi x_{i}) + \varepsilon_{i} \, ,
\end{align}
which are all distorted by artificial noise $\varepsilon_{i}$. At each ordinate $y_{i}$, the
latter is a random variable drawn from a normal distribution centered about zero with variance
$(0.03)^2$, i.e., $\varepsilon_{i} \sim \calN(\mu = 0, \sigma^{2} = (0.03)^{2})$.

In order to compare results, we use either the absolute percentage error (APE)
\begin{equation}
	\label{eq:sample_data_03}
	\ape(a,b) = \frac{\lvert a - b \rvert}{\max\{ \lvert a \rvert, \delta \}} \times 100 \, ,
\end{equation}
where $a,b \in \setR$ and $\delta > 0$, or the root mean squared error (RMSE)
\begin{equation}
	\label{eq:sample_data_04}
	\rmse(f^{\opt}) = \sqrt{\frac{1}{n} \sum_{i=0}^{n-1} (y_{i} - f^{\opt}(x_{i}))^2} \, ,
\end{equation}
where $f^{\opt}(x) = \sum_{j=0}^{m-1} c_{j}^{\opt} \phi_{j}(x)$ is a given optimal fit with
optimal expansion coefficients $c_{j}^{\opt}$ obtained from solving the least-squares optimization
problem (see Eqs.~\labelcref{eq:least_squares} and \labelcref{eq:function_expansion}). Equation
\labelcref{eq:sample_data_04} aggregates the squared differences between the actual and
approximated values for each sample data point---the very same measure that is minimized by the
method of least squares---and is thus well-suited to compare optimal least-squares fits $f^{\opt}$
obtained via different techniques (see Sections \ref{curve_fitting_classical} and
\ref{curve_fitting_qubo}).

Furthermore, we always normalize the ordinate values to the unit interval before the optimal
fit is computed. That is, the least-squares method actually operates on min-max-normalized data
points $\{ (x_{i}, y_{i}') \}_{i=0}^{n-1}$ with
\begin{equation}
	\label{eq:sample_data_05}
	y_{i}' = \frac{y_{i} - y_{\textup{min}}}{y_{\textup{max}} - y_{\textup{min}}} \in [0, 1] \, ,
\end{equation}
where $y_{\textup{min}} = \min_{i} y_{i}$ and $y_{\textup{max}} = \max_{i} y_{i}$.
Thus, all optimal expansion coefficients $c_{j}^{\opt}$ presented here are understood to be obtained
with respect to min-max-normalized data $y_{i}'$.

Now, in the following we consider three different standardized functions $\phi_{j}$ for the
expansion \labelcref{eq:function_expansion} of the fit function, apply them to our sample data,
and discuss the quality of the solutions obtained via different strategies as well as the influence
of certain parameters.

\subsection{Triangular functions}
\label{results_triangular}

\begin{table*}[t]
	\centering%
	\caption{\label{tab:linear_regression}%
		Optimal linear regression parameters obtained classically and by solving a QUBO.
		Also shown are the absolute percentage errors (APE) of the optimal
		parameters obtained from the QUBO with respect to the classical ones.
	}%
	\vspace*{0.25em}%
	\begin{ruledtabular}%
	\begin{tabular}{lcccc}%
		& $c_{0}^{\opt}$ & $c_{1}^{\opt}$ & APE $c_{0}^{\ast}$ [\%] & APE $c_{1}^{\ast}$ [\%] \\%
		\cmidrule{2-5}%
		Classical (matrix inversion) & $0.010664$ & $0.978332$ & -- & -- \\%
		QUBO (tabu search) & $0.010742$ & $0.978516$ & $0.73$ & $0.02$ \\%
		QUBO (quantum annealing) & $0.011719$ & $0.977539$ & $9.89$ & $0.08$ \\%
	\end{tabular}%
	\end{ruledtabular}%
\end{table*}

\begin{table*}[t]
	\centering%
	\caption{\label{tab:linear_regression_2}%
		Same as \cref{tab:linear_regression} but with the linear sample data now given by
		$y_{i} = -x_{i}/4 + 1/3 + \varepsilon_{i}$.
	}%
	\vspace*{0.25em}%
	\begin{ruledtabular}%
	\begin{tabular}{lcccc}%
		& $c_{0}^{\opt}$ & $c_{1}^{\opt}$ & APE $c_{0}^{\ast}$ [\%] & APE $c_{1}^{\ast}$ [\%] \\%
		\cmidrule{2-5}%
		Classical (matrix inversion) & $0.874487$ & $0.200203$ & -- & -- \\%
		QUBO (tabu search) & $0.875000$ & $0.200195$ & $0.06$ & $0.004$ \\%
		QUBO (quantum annealing) & $0.872070$ & $0.196289$ & $0.28$ & $1.96$ \\%
	\end{tabular}%
	\end{ruledtabular}%
	\vspace*{1em}%
	\caption{\label{tab:n_dependence}%
		APE between the RMSE of the optimal classical solution $\foptcl$ and the optimal QUBO
		solutions $\fopttb$ and $\foptqa$ obtained via tabu search and quantum annealing, respectively,
		for different values of $n$. The underlying sample data set is the linear one
		(\cref{eq:sample_data_02a}) and a linear regression is used for the approximation.
	}%
	\vspace*{0.25em}%
	\begin{ruledtabular}%
	\begin{tabular}{ccc}%
		$n$ & $\ape(\rmse(\foptcl), \rmse(\fopttb))$ [\%] & $\ape(\rmse(\foptcl), \rmse(\foptqa))$ [\%] \\%
		\cmidrule{1-3}%
		$64$ & $2.88 \times 10^{-4}$ & $2.23 \times 10^{-2}$ \\%
		$128$ & $3.64 \times 10^{-4}$ & $4.50 \times 10^{-2}$ \\%
		$256$ & $1.10 \times 10^{-4}$ & $1.22 \times 10^{-2}$ \\%
		$512$ & $1.09 \times 10^{-4}$ & $1.71 \times 10^{-2}$ \\%
		$1024$ & $6.80 \times 10^{-4}$ & $3.16 \times 10^{-2}$ \\%
	\end{tabular}%
	\end{ruledtabular}%
\end{table*}

We begin with triangular functions as defined in \cref{eq:triangular}, $\phi_{j} = \triang_{j}$,
i.e., the fit function is given by
\begin{equation}
	\label{eq:results_triang_01}
	f(x) = c_{0} \triang_{0}(x) + \dotsb + c_{m-1} \triang_{m-1}(x) \, ,
\end{equation}
and the least squares optimization problem thus explicitly reads
\begin{equation}
	\label{eq:results_triang_02}
	\min_{c \in \setR^{m}}
	\sum_{i=0}^{n-1} \left( y_{i} - \sum_{j=0}^{m-1} c_{j} \triang_{j}(x_{i}) \right)^{2} \, .
\end{equation}

First, we consider the simple yet important case of linear regression. To this end, we choose
$n = 64$ and the linear sample data \labelcref{eq:sample_data_02a}. Furthermore, we set $m = 2$,
which implies that the two supporting points of the triangular functions (see \cref{eq:triangular})
are simply the first and last abscissa values \labelcref{eq:sample_data_01} of the sample data,
respectively, i.e., $\tilde{x}_{0} = x_{0} = 0$ and $\tilde{x}_{1} = x_{n-1} = 1$. This results in
a fit function $f(x) = c_{0} \triang_{0}(x) + c_{1} \triang_{1}(x)$ that is a straight line across
the whole abscissa interval $[0, 1]$ of our sample data. For the corresponding QUBO formulation of
the problem, we choose $d = 10$ and $p = d - 2 = 8$, which results in a QUBO matrix of size
$md \times md = 20 \times 20$.

In \cref{tab:linear_regression}, we show our results for the optimal linear regression parameters
$c_{0}^{\opt}$ and $c_{1}^{\opt}$ obtained classically by (naive) matrix inversion according to
\cref{eq:classical_04} and by solving the corresponding QUBO as given in \cref{eq:quantum_06};
we also show the APE (see \cref{eq:sample_data_03}) of the optimal coefficients obtain via a QUBO
with respect to the classical results. The QUBO is solved by means of a tabu search as well as by quantum
annealing on a D-Wave quantum processing unit (QPU).%
\footnote{For the tabu search, we used the \textit{MST2} multi-start tabu algorithm of
Ref.~\cite{Palubeckis:2004a}. Regarding quantum annealing for solving the QUBO, we used
D-Wave's \textit{Advantage2} system. In particular, we performed 100 runs with 100 reads each with
an annealing time of 200 microseconds to collect potential solution candidates, yielding 10000 samples
in total, and finally took the lowest-energy solution. This setup is used throughout the whole section
whenever a QUBO problem is solved via tabu search or quantum annealing.}
Comparing the three solutions, we find that all agree up to the second decimal point. The QUBO
solution found by the tabu search is particularly close to the exact classical solution, agreeing
up to the third decimal point. Though the quantum annealing solution differs more noticeable, it
is still in good agreement with the two other results. Furthermore, we find qualitatively similar
results if we replace the numerical coefficients in the linear sample data \labelcref{eq:sample_data_02a}
by other values, e.g., randomly drawn from $[-1, 1]$. As an example, \cref{tab:linear_regression_2}
shows the results for a different set of linear sample data.

Thus, in the case of linear regression, our results indicate that formulating curve fitting as a
QUBO problem indeed works and yields a solution in satisfying agreement with the classical one.
They also suggest that quantum annealing in its current state is in principle able (within limits)
to deliver promising and practically usable results for this continuous problem.%
\footnote{Because of the probabilistic (or stochastic) nature of quantum annealing and the fact
that the QPU is a real (i.e., non-ideal) system---thus susceptible to thermal fluctuations and
interactions with its environment---subsequent calculations with the same input may yield slightly
different results. Based on our experiences and observations in the course of this work, the
quantum annealing results presented here seem to be accompanied with an overall error of about
five to ten percent.}

It is furthermore instructive to investigate the influence of certain parameters on the solutions.
In \cref{tab:n_dependence}, we show the APE between the RMSE of the optimal classical
solution $\foptcl$ and the optimal QUBO solutions $\fopttb$ and $\foptqa$ obtained via tabu search
and quantum annealing, respectively, for different values of the number of sample data points $n$
(again for the case of linear regression with sample data \labelcref{eq:sample_data_02a}). We
conclude that the number of to-be-approximated data points has a negligible influence on the
solution quality. This is expected because $n$ enters only implicitly in the course of determining
the entries of the QUBO matrix (see Eqs.~\labelcref{eq:classical_01a} and \labelcref{eq:classical_01b})
while its size is given by $md \times md$, see \cref{eq:quantum_07}, i.e., only determined by
$m$ and $d$. We find the same result if we use our other sample data sets instead of the linear one.
Therefore, we fix $n = 64$ for the remainder of this section.

\begin{figure}[t]
	\centering%
	\includegraphics[scale=0.9]{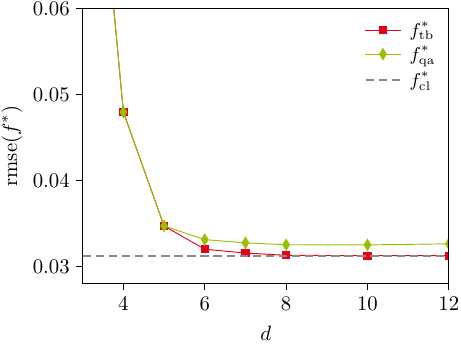}%
	\vspace*{-0.5em}%
	\caption{\label{fig:triangular_d_dependence}%
		RMSE of the tabu search (red squares, $\fopttb$) and quantum annealing (green diamonds, $\foptqa$)
		solutions as a function of $d$. The RMSE of the exact, $d$-independent classical solution $\foptcl$
		is shown as a gray dashed line. Lines between data points are shown to guide the eye. The underlying
		sample data set is the cubic one (\cref{eq:sample_data_02c}), and triangular functions with
		$m=4$ are used for the approximation.
	}%
\end{figure}

In \cref{fig:triangular_d_dependence}, we show the RMSE of the tabu (red squares) and quantum
annealing (green diamonds) solutions $\fopttb$ and $\foptqa$, respectively, as a function of
$d$, the number of digits that is used for the binary fixed-point format of the coefficients $c_{j}$,
see \cref{eq:quantum_02a}, with $p = d - 1$ fixed. Furthermore, we show the RMSE of the classical
solution $\foptcl$, which is exact in the sense of numerical precision and thus independent of $d$.
Here, we use the cubic sample data \labelcref{eq:sample_data_02c} and $m = 4$, i.e., we are beyond
linear regression and have a truly piecewise-linear fit of the form
$f(x) = c_{0} \triang_{0}(x) + \dotsb + c_{3} \triang_{3}(x)$ that is not just a straight regression
line. For $d \geq 8$, the tabu solutions are very close to the classical one and a further
increase of $d$ has only a negligible impact. The same holds for the quantum annealing solutions,
though they are always less accurate than the corresponding tabu results. Both are only identical
provided the QUBO matrix is small enough, namely for $d \leq 5$, which corresponds to a QUBO matrix
of size $20 \times 20$ or smaller. In our opinion, however, the quantum annealing solutions for
larger QUBO sizes are still in acceptable agreement with the tabu and classical ones. We find
qualitatively similar results if we use the quadratic sample data or if the numerical coefficients
of the cubic data are randomly drawn from $[-1, 1]$. Therefore, $d = 8$ (with $p = d - 1$)
seems to be sufficient for the binary fixed-point representation of the expansion coefficients $c_{j}$
regarding the cases considered in this work. Furthermore, $d \in \{ 6, 7 \}$ are reasonable compromises
between accuracy and keeping the QUBO matrix as small as possible.

\begin{figure*}[p]
	\centering%
	\includegraphics[scale=0.9]{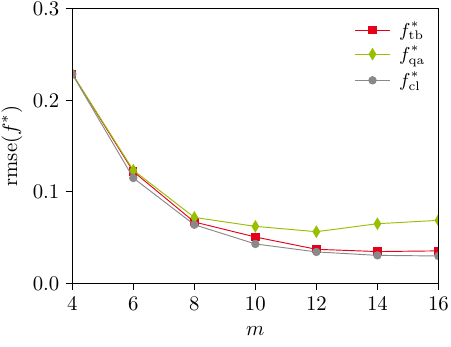}%
	\hspace*{5pt}%
	\includegraphics[scale=0.9]{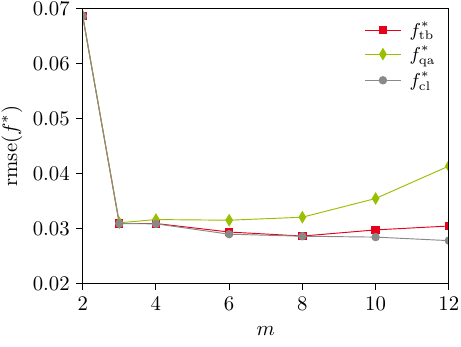}%
	\vspace*{-0.5em}%
	\caption{\label{fig:triangular_m_dependence}%
		RMSE of the tabu search (red squares, $\fopttb$), quantum annealing (green diamonds, $\foptqa$)
		and classical (gray circles, $\foptcl$) solutions as a function of $m$. The underlying sample
		data are the trigonometric (left diagram) and quadratic (right diagram) ones, respectively.
	}%
\end{figure*}

Last, we consider the case of varying $m$, i.e., how the solution quality depends on the number
of triangular functions used for the curve fitting. This is an particularly interesting case
because one usually chooses a sufficient numerical accuracy and simply increases $m$ to the point
where the data is sufficiently good approximated. Curve fitting with triangular functions is well-suited
for that because the result is a piecewise-linear approximation (a polygonal chain), which is often
used to approximate more complex curves or if the functional form of the to-be-approximated data
is not immediately obvious. In \cref{fig:triangular_m_dependence}, we show the RMSE of the tabu
(red squares), quantum annealing (green diamonds), and classical (gray circles) solutions as a function
of $m$, where the sample data sets are the trigonometric (left diagram) and quadratic ones (right diagram);
again, $n = 64$, $d = 8$, and $p = d - 1 = 7$. First, we consider the left diagram. It shows the overall
trend that the RMSE of the classical solution becomes smaller as $m$ increases. Initially, the tabu
and quantum annealing solutions follow that trend. For $m \geq 10$ (which corresponds to a QUBO matrix
of size $80 \times 80$), however, the quantum annealing solution becomes worse with increasing $m$.
The right diagram, though with a rather different underlying set of sample data (quadratic instead
of trigonometric), shows the same qualitative behavior. Quantitatively, noticeable differences and
a continuous worsening of the quantum annealing solutions with increasing $m$ already starts at
$m = 8$, which corresponds to a QUBO matrix size of $64 \times 64$.

We also tested other polynomial sample data as well as varied the numerical coefficients of the sample data
used for \cref{fig:triangular_m_dependence} and found qualitatively similar results. This indicates
that a QUBO size of about $60 \times 60$ to $80 \times 80$ seems to be the range above which D-Wave's
current quantum annealing hardware fails to find a satisfying solution of the QUBO formulation of least
squares curve fitting with triangular functions.

\subsection{Chebyshev polynomials}
\label{results_chebyshev}

Besides the triangular functions, we shall also consider Chebyshev polynomials of the first kind as
defined in \cref{eq:chebyshev}, $\phi_{j} = T_{j}$, i.e., the fit function is now given by
\begin{equation}
	f(x) = c_{0} + \sum_{j=1}^{m-1} c_{j} T_{j}(x)
\end{equation}
and is a polynomial of degree $m$.

\begin{table*}[p]
	\centering%
	\caption{\label{tab:chebyshev_quad}%
		Optimal expansion coefficients of the individual Chebyshev polynomials of a quadratic fit,
		obtained classically via matrix inversion and by solving a QUBO via tabu search as well as
		quantum annealing (labeled ``qa''). Also shown are the APE of the optimal coefficients obtained
		from the QUBO with respect to the classical ones.
	}%
	\vspace*{0.25em}%
	\begin{ruledtabular}%
	\begin{tabular}{lcccccc}%
		& $c_{0}^{\opt}$ & $c_{1}^{\opt}$ & $c_{2}^{\opt}$ & APE $c_{0}^{\ast}$ [\%] & APE $c_{1}^{\ast}$ [\%] & APE $c_{2}^{\ast}$ [\%] \\%
		\cmidrule{2-7}%
		Classical & $0.544598$ & $-0.031799$ & $0.499257$ & -- & -- & -- \\%
		QUBO (tabu) & $0.523438$ & $-0.039063$ & $0.484375$ & $3.89$ & $22.84$ & $2.98$ \\%
		QUBO (qa) & $0.445313$ & $-0.101563$ & $0.421875$ & $18.23$ & $219.39$ & $15.50$ \\%
	\end{tabular}%
	\end{ruledtabular}%
	\vspace*{2.5em}%
\end{table*}

\begin{figure*}[p]
	\centering%
	\includegraphics[scale=0.9]{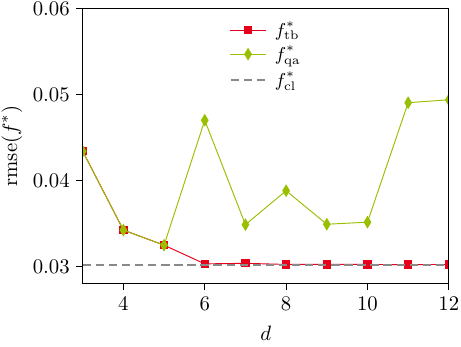}%
	\vspace*{-0.5em}%
	\caption{\label{fig:chebyshev_d_dependence}%
		RMSE of the tabu search (red squares, $\fopttb$) and quantum annealing (green diamonds, $\foptqa$)
		solutions as a function of $d$. The RMSE of the exact, $d$-independent classical solution $\foptcl$
		is shown as a gray dashed line. Lines between data points are shown to guide the eye. The underlying
		sample data set is the cubic one (\cref{eq:sample_data_02c}), and Chebyshev polynomials with
		$m=4$ are used for the approximation.
	}%
\end{figure*}

\begin{figure*}[t]
	\centering%
	\includegraphics[scale=0.9]{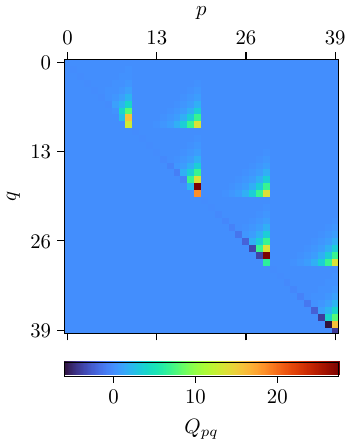}%
	\hspace*{1em}%
	\includegraphics[scale=0.9]{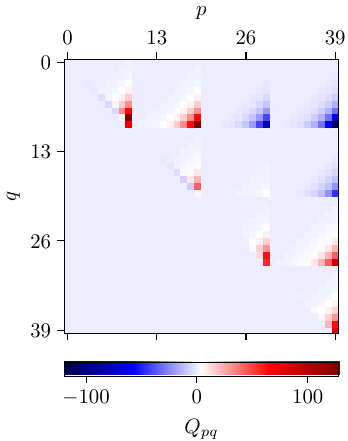}%
	\vspace*{-0.5em}%
	\caption{\label{fig:qubo_matrices}%
		Heat map plots of exemplary QUBO matrices $Q$ of piecewise-linear approximation
		with triangular functions (left) and Chebyshev polynomial approximation (right).
		Here, we display the upper-triangularized versions of these matrices since every QUBO
		problem can be written in a way such that $Q$ is an upper-triangular matrix.
	}%
\end{figure*}

In \cref{tab:chebyshev_quad}, we show our results for fitting a quadratic Chebyshev polynomial
(i.e., $m = 3$) to the quadratic sample data \labelcref{eq:sample_data_02b}; again, $n = 64$, $d = 8$,
and $p = d - 1 = 7$, which results in a QUBO matrix of size $24 \times 24$. While the tabu solution
is in acceptable agreement with the classical one, the quantum annealing solution deviates significantly.
We find the same result for the cubic data (also with varying numerical coefficients) with a cubic Chebyshev
polynomial ($m = 4$) as the fit function. In \cref{fig:chebyshev_d_dependence}, we show the RMSE
of the tabu (red squares) and quantum annealing (green diamonds) solutions $\fopttb$ and $\foptqa$,
respectively, as a function of $d$, the number of digits that is used for the binary fixed-point format
of the expansion coefficients. Furthermore, we show the RMSE of the classical solution, which is
independent of $d$, as a gray dashed line. The underlying sample data is the cubic one and the Chebyshev
fit is of degree three, i.e., $m=4$. We plot against $d$ because we would like to keep the degree of the
Chebyshev polynomial fixed while letting the QUBO matrix increase in size. While the tabu solution
shows the expected behavior (cf.~\cref{fig:triangular_d_dependence}), we are not able to obtain
a satisfying quantum annealing solution for problem sizes larger than $20 \times 20$. This is in stark
contrast to our findings for piecewise-linear curve fitting using triangular functions.

This seems to hint toward the limitations of D-Wave's current quantum annealing hardware for
curve fitting with polynomials formulated as a QUBO. More precisely, the corresponding QUBO problem,
which can be seen as a graph with logical qubits representing nodes that are connected by edges with
weights given by the matrix elements of the QUBO matrix, needs to be mapped onto the so-called working
graph of the quantum annealing device. Therefore, one might need several physical qubits to represent one
logical qubit. Particularly for highly-connected (dense) QUBO problems, this could yield an Ising model
that exhibits a complicated low-energy landscape where the ground state is difficult to access. This would
explain why we were not successful in obtaining a feasible Chebyshev approximation as soon as
its QUBO matrix is larger than the above mentioned size of about $20 \times 20$. Such a matrix is,
as apparent from the right diagram of \cref{fig:qubo_matrices}, indeed relatively highly connected.

The amount of nonvanishing off-diagonal elements of the QUBO matrix for least-squares curve
fitting is solely determined by the standardized functions $\phi_{j}$ (keeping the underlying
to-be-approximated data fixed); see \cref{eq:classical_01b}. More precisely, in case of the
Chebyshev approximation we deal with polynomials $T_{j}$ whose support are the whole real line,
\begin{equation}
	\label{eq:numerics_04}
	\supp(T_{j}) = \setR
\end{equation}
for all $j$. Therefore, the resulting QUBO problem corresponds to a
rather highly connected graph, and the matrix has many nonzero elements even far off the diagonal;
see the right diagram of \cref{fig:qubo_matrices} for a graphical depiction. On the other hand,
the QUBO matrix for a piecewise-linear approximation is always a band matrix, see the left diagram
of \cref{fig:qubo_matrices}, because the triangular functions $\triang_{j}$ have compact support,
\begin{equation}
	\label{eq:numerics_05}
	\supp(\triang_{j})
	=
	\begin{cases}
		[\tilde{x}_{0}, \tilde{x}_{1}] & \text{if~} j = 0 \, ,
		\\[0.25em]
		[\tilde{x}_{m-2}, \tilde{x}_{m-1}] & \text{if~} j = m - 1 \, ,
		\\[0.25em]
		[\tilde{x}_{j-1}, \tilde{x}_{j+1}] & \text{otherwise} \, .
	\end{cases}
\end{equation}
In our opinion, this explains why
\begin{itemize}[leftmargin=*,widest=(ii)]
	\item[(i)] we are unsuccessful in finding a satisfying quantum annealing solution
	for the Chebyshev approximation for QUBO matrices larger than a certain size;
	\item[(ii)] the Chebyshev solution is less accurate compared to piecewise-linear
	approximation for feasible cases of the size of the QUBO problem.
\end{itemize}

Furthermore, another obstacle might be the order of magnitude of the elements of the QUBO
matrix. These have to be mapped onto the coupling strengths and qubit biases of the Ising model,
which are constrained by the physical limits of the quantum annealing hardware, i.e., the
available ranges of couplings and biases as well as their coarseness. As a consequence, it might be
the case that a complicated low-energy landscape cannot be ``scanned'' accurate enough by the
current quantum annealing devices. Unfortunately, the actual values of these matrix elements
are determined by the to-be-approximated data and the employed standardized functions $\phi_{j}$
and thus cannot be freely adjusted or altered.

We would like to emphasize that by solving the QUBO problems for piecewise-linear
and Chebyshev approximations using a heuristic technique, e.g., tabu search or simulated
annealing, we (almost) always find a solution in satisfying agreement with the classical one---even
for sizes of the QUBO matrix for which the quantum annealing solution turned out to be unpractical.
Thus, we are confident that a next-generation quantum annealing device will be capable of solving
such QUBO problems.

Finally, we would like to mention that we always subjected the ``raw'' QUBO matrix to the quantum
annealing process. It is well-known that dense QUBO problems are hard to solve, and a proper
preprocessing of the QUBO matrix might improve the solution quality. Possible preprocessing
steps are sparsification of the QUBO matrix or properly modifying/combining matrix elements
\cite{Lewis:2017,Glover:2017,Suppakitpaisarn:2024}, possibly resulting in a reduction in size
and/or density of the problem. Another approach would be quantum annealing-inspired algorithms
like SimCIM \cite{Tiunov:2019} or the NGQ solver \cite{Du:2025}, which seem to be more robust
against dense QUBO problems \cite{Hamerly:2018,Takesue:2025,Du:2025}.
However, these investigations are beyond the scope of this work and deferred to possible future work.

\section{Speed profile optimization with dynamic programming}
\label{speed_profile}

In this section, we touch upon a real-world use case, with which we are
currently concerned, where curve fitting can come in handy. Namely, a cargo vessel's
voyage under the constraint of minimal fuel consumption while arriving at the destination at
the requested time of arrival. Phrased differently, the aim is to find an optimized speed
profile that minimizes fuel costs while arriving just in time. From our point of view, this task
lends itself to be solved within the framework of dynamic programming
\cite{Bellman:1952a,Bellman:1953a,Bellman:1954a}; see, e.g., Refs.~\cite{Bellman:1957a,%
Neumann:2002a,Bertsekas:2017a,Bertsekas:2012a,Sutton:2018a,Bertsekas:2019a} for comprehensive
overviews. Work in that direction has been done regarding the optimal control of sailboats
\cite{Ferretti:2017a,Miles:2021a,Wang:2023a}. To our knowledge, however, dynamic programming is
not used for maritime just-in-time navigation of cargo vessels yet, and we would like to leverage
its use.

\subsection{General formalism}
\label{speed_profile_general}

More formally speaking, we consider a deterministic Markov decision process that starts at
time $t = 0$ and terminates at $t = T$, where the time evolution is discrete in terms of integer
time steps: $t = 0, \dotsc, T$. The process is characterized by sets $\setStates{t} \subset \setR^{k}$
and $\setActions{t} \subset \setR^{l}$, where $k,l \in \setN$, that contain all states the process can
be in and all actions that can be performed in order to change the current state, respectively,
at time step $t$.%
\footnote{In general, the actions depend on the state $\mystate{t} \in \setStates{t}$ of the 
process, i.e., $\setActions{t} = \setActions{t}(\mystate{t})$. However, in terms of notation, we omit this
dependency for the the sake of brevity.}

The dynamics of the decision process, i.e., how the next state $\mystate{t+1} \in \setStates{t+1}$
is obtained from the previous one $\mystate{t} \in \setStates{t}$ by performing a certain action
$\action{t} \in \setActions{t}$, are governed by a transition function
$\transFunc{t} \colon \setStates{t} \times \setActions{t} \to \setStates{t+1}$; that is,
\begin{equation}
	\label{eq:speed_profile_01}
	\mystate{t+1} = \transFunc{t}(\mystate{t}, \action{t}) \, .
\end{equation}

The immediate costs that occur if the action $\action{t}$ is performed while being in state $\mystate{t}$
are given by a cost function $\costFunc{t} \colon \setStates{t} \times \setActions{t} \to \setR$.
This motivates the definition of a value function $\valFunc{t} \colon \setStates{t} \times
\setActions{t} \times \dotso \times \setActions{T-1} \to \setR$, which is given by
\begin{equation}
	\label{eq:speed_profile_02}
	\valFunc{t}(\mystate{t}, \action{t}, \dotsc, \action{T-1})
	= \sum_{s=t}^{T-1} \costFunc{s}(\mystate{s}, \action{s}) + D(\mystate{T}) \, ,
\end{equation}
where $D \colon \setStates{T} \to \setR$ encodes terminal costs that occur at the end of the
process. The value function is a measure for how cost-effectively the decision process has been
if started at time step $t$ with initial state $\mystate{t}$ and running until the end by virtue of a
given sequence of actions $\action{t}, \dotsc, \action{T-1}$. The states $\mystate{t+1}, \dotsc, \mystate{T}$
in \cref{eq:speed_profile_02} are generated according to the temporal evolution defined by the
process' dynamics \labelcref{eq:speed_profile_01}.

The goal is to find optimal actions $\action{t}^{\opt}, \dotsc, \action{T-1}^{\opt}$ that minimize
the cost function for a given initial state. Consequently, these optimal actions, which can be
regarded as a recipe how to operate the process such that minimal total costs incur, are encoded
in the optimal value function $\valFuncOpt{t} \colon \setStates{t} \to \setR$, which reads
\begin{equation}
	\label{eq:speed_profile_03}
	\valFuncOpt{t}(\mystate{t}) =
	\min_{\substack{\action{s} \in \setActions{s} \\ s = t, \dotsc, T - 1}}
	\valFunc{t}(\mystate{t}, \action{s}, \dotsc, \action{T-1}) \, .
\end{equation}
For any given initial state, it yields the optimal remaining costs; usually, one is interested
in $\valFuncOpt{0}$, i.e., the total optimal costs for the whole process. Now, it can be shown that
$\valFuncOpt{t}$ obeys the so-called Bellman equation \cite{Bellman:1952a,Bellman:1953a,%
Bellman:1954a,Bellman:1957a}, which lies at the heart of the dynamic programming paradigm and plays
an important role in the field of reinforcement learning \cite{Sutton:2018a,Bertsekas:2019a}.
It relates the optimal value function at different time steps in a recursive manner and is given by
\begin{equation}
	\label{eq:speed_profile_04}
	\valFuncOpt{t}(\mystate{t}) = \min_{\action{t} \in \setActions{t}}
	\bigl\{ C_{t}(\mystate{t}, \action{t}) + \valFuncOpt{t+1}(\mystate{t+1}) \bigr\}
\end{equation}
for $t = 0, \dotsc, T - 1$ and all $\mystate{t} \in \setStates{t}$ with
$\mystate{t+1} = \transFunc{t}(\mystate{t}, \action{t})$ and boundary condition
$\valFuncOpt{T}(\mystate{T}) = D(\mystate{T})$.

Solving the recursively-coupled set of equations is a nontrivial task. Fortunately, the Bellman
equation suggests a practical algorithm to obtain the globally optimal solution. Typically, one
is interested in the optimal control $\action{0}^{\opt}, \dotsc, \action{T-1}^{\opt}$ of the whole
process. The solution strategy consists of a backward (value) and forward (policy) iteration;
see, e.g., Ref.~\cite{Neumann:2002a}. First, one computes the function values of the optimal value function
starting from the last time step and then moving sequentially backward in time, i.e., for all
$t = T, \dotsc, 0$ and for all states $\mystate{t} \in \setStates{t}$, the values
$\valFuncOpt{t}(\mystate{t})$ are computed according to \cref{eq:speed_profile_04}. Second, the
optimal actions are obtained via
\begin{equation}
	\label{eq:speed_profile_05}
	\action{t}^{\opt} = \argmin_{\action{t} \in \setActions{t}}
	\bigl\{ \costFunc{t}(\mystate{t}^{\opt}, \action{t})
	+ \valFuncOpt{t+1}(\transFunc{t}(\mystate{t}^{\opt}, \action{t})) \bigr\}
\end{equation}
with $\mystate{t+1}^{\opt} = \transFunc{t}(\mystate{t}^{\opt}, \action{t}^{\opt})$ and initial state
$\mystate{0} = \mystate{0}^{\opt}$ for $t = 0, \dotsc, T - 1$, i.e., forward in time.

\subsection{Toward an optimized speed profile}
\label{speed_profile_specific}

In principle, the Bellman equation as given in \cref{eq:speed_profile_04} can be used right away
to compute an optimized speed profile. The mapping onto that task is, for example, accomplished by
the following setup:
\begin{itemize}[leftmargin=*,label=\textbullet]
	\item The trajectory of the vessel is fixed and covers a certain distance of $\ell$
	nautical miles.
	\item $\setStates{t}$ contains all feasible positions of the vessel on the trajectory
	(i.e., the distances traveled since departure at time $t = 0$) and environmental aspects
	like weather, currents, and forecasts at time $t$.
	\item $\setActions{t}$ contains all feasible speeds (through water) of the vessel at time $t$.
	\item $\costFunc{t}(\mystate{t}, \action{t})$ are the fuel costs if the vessel moves
	with speed $\action{t}$ while being at position $\mystate{t}$ at time $t$.
	\item $\transFunc{t}(\mystate{t}, \action{t})$ describes how a speed of $\action{t}$ while being
	at position $\mystate{t}$ at time $t$ translates into a change in position to $\mystate{t+1}$.
\end{itemize}
Given that all of the above are available, particularly the cost and transition functions for
arbitrary $\mystate{t}$ and $\action{t}$, solving the Bellman equation yields optimal speeds and
positions $\action{0}^{\opt}, \dotsc, \action{T-1}^{\opt}$ and $\mystate{0}^{\opt}, \dotsc, \mystate{T}^{\opt}$,
respectively, on the trajectory that globally minimize the value function, i.e., provide the
optimal operating scenario for just-in-time arrival at minimal fuel consumption.

\subsection{Proof of principle}
\label{speed_profile_example}

For the purpose of a proof of principle, we choose a linear-quadratic optimal control problem. The
system's dynamics are linear in both the state and action variables, while the cost functional comprises
both running and terminal costs, each formulated as quadratic forms. This specific structure offers several
analytical and computational advantages. First, the convexity of the cost functional and the linearity
of the system's dynamics ensure that the resulting optimization problem is convex. This implies the
existence of an unique global optimum. Second, the problem admits an analytical solution via the Riccati
differential equation, enabling the derivation of the optimal control law in the form of a linear state
feedback.

In the following, we investigate a simple implementation of what we described above in the previous
subsection, namely the problem of just-in-time arrival. We set
\begin{align}
	\label{eq:speed_profile_06}
	\setStates{t}, \setActions{t} &\subset \setR \, ,
	\\[0.25em]
	\label{eq:speed_profile_07}
	\costFunc{t}(\mystate{t}, \action{t}) &= \biggl( \frac{\action{t} - w(\mystate{t})}{\vmax} \biggr)^{2} \, ,
	\\[0.25em]
	\label{eq:speed_profile_08}
	\mystate{t+1} &= \transFunc{t}(\mystate{t}, \action{t}) = \mystate{t} + \action{t} - w(\mystate{t}) \, ,
	\\[0.25em]
	\label{eq:speed_profile_09}
	D(\mystate{T}) &= \alpha \, \Bigl( 1 - \frac{\mystate{T}}{\ell} \Bigr)^{2} + 1
\end{align}
for all $t$.%
\footnote{From a physics point of view, $\action{t}$ and $w$ are velocities. Thus, the transition
function actually reads $\mystate{t+1} = \mystate{t} + (\action{t} - w(\mystate{t})) \+ \Delta t$
with $\Delta t$ denoting the time difference between two time steps. Since we are working solely with
dimensionless quantities and integer time steps, $\Delta t = 1$ and it is omitted for the sake of brevity.}
Here, $w \colon \setStates{t} \to \setR$ denotes the magnitude of the projected sea current
with respect to the direction of travel at position $\mystate{t}$ on the trajectory, $\vmax$ represents
the maximal possible speed through water, and $\alpha \in \setRppos$ is a numerical parameter to
fine-tune the terminal costs.

Upon inserting these definitions into \cref{eq:speed_profile_02} and writing the final state as
$\mystate{T} = \mystate{s} + \sum_{t=s}^{T-1} (\action{t} - w(\mystate{t}))$ as a sum of actions
starting at time $s$ in state $\mystate{s}$, which follows directly from
\cref{eq:speed_profile_08}, the value function for the whole process reads
\begin{align}
	\label{eq:speed_profile_10}
	\valFunc{s}(\mystate{s}, \action{s}, \dotsc, &\,\action{T-1})
	=
	\frac{1}{\vmax^{2}} \sum_{t=s}^{T-1} (\action{t} - w(\mystate{t}))^{2}
	\notag\\
	& +\alpha
	\biggl[ \ell - \frac{1}{\ell}  \biggl( \mystate{s} - \sum_{t=s}^{T-1} (\action{t} - w(\mystate{t})) \biggr) \biggr]^{2}
	\notag\\
	& + 1 \, .
\end{align}

Without loss of generality (and for numerical convenience), we assume that the sets of states and
actions $\setStates{t}$ and $\setActions{t}$, respectively, are the same for every time step, i.e.,
independent of $t$. It is then straightforward to find the optimal actions, which are the roots of
$\grad_{\? u} \valFunc{s}$ (for fixed $\mystate{s}$). Regarding the gradient of the value function,
we find
\begin{equation}
	\label{eq:speed_profile_10_1}
	\frac{\partial \valFunc{0}}{\partial \action{t}}
	=
	\frac{2\action{t}}{\vmax^{2}}
	-
	\frac{2\alpha}{\ell^{2}}
	\biggl[ \ell - \mystate{s} - \sum_{s^{\prime}=s}^{T-1} (\action{s^{\prime}} - w(\mystate{s^{\prime}})) \biggr] \, .
\end{equation}

Now, we shall solve this optimization problem by numerically solving the Bellman equation
\labelcref{eq:speed_profile_04} and subsequently using \cref{eq:speed_profile_05} to find
the optimal actions. For the sake of simplicity, we set $w(\mystate{t}) = 0$ for all
$\mystate{t} \in \setStates{t}$ and all $t$, and consider a voyage with $\ell = 100$, $\vmax = 50$,
and $T = 4$. Neglecting phases of (de)acceleration and without currents, the optimal solution
is obviously $\action{0}^{\opt} = \dotso = \action{3}^{\opt} = 25$, i.e., moving (trivially)
with constant velocity. Approaching this problem in terms of the value function
\labelcref{eq:speed_profile_09}, we find from \cref{eq:speed_profile_10_1} the analytical solution
\begin{equation}
	\label{eq:speed_profile_11}
	\action{t}^{\opt} = \frac{\alpha \vmax^{2} (\ell - \mystate{s})}{T - s + \ell^{2}}
\end{equation}
for all $t = s, \dotsc, T$. As expected, this analytical solution describes a constant
time-independent velocity, which consistently coincides with the optimal solution in the limit
$\alpha \to \infty$, i.e., for infinite terminal costs. This analytical solution (with finite $\alpha$)
shall serve for comparison with results from solving the Bellman equation numerically.

In the course of that, the traditional/naive approach discretizes the state space. Hence, after
the value iteration, the resulting optimal value function
$\valFuncOpt{t}$ is defined only on that discrete state space grid. Consequently, this poses
constraints on the possible discretizations of the action space for the policy iteration
\labelcref{eq:speed_profile_05} because for a given $\mystate{t}$, $\action{t}$ must be chosen
such that $\mystate{t+1} = \transFunc{t}(\mystate{t}, \action{t})$ lies on the state space grid,
too, or the off-grid result $\transFunc{t}(\mystate{t}, \action{t})$ has to be ``snapped'' onto
the nearest grid point. Both approaches induce a discretization error that manifests in the
final (pseudo-)optimal policy. The in that way obtained optimal policy most likely differs
from the true, globally-optimal solution---it is only optimal with respect to the discretization
of the state space. In our case, this would result in a pseudo-optimal speed profile different
from the analytical solution \labelcref{eq:speed_profile_11}.

In principle, this issue can be mitigated by using a sufficiently dense state space grid.
In our opinion, however, this is not feasible for many real-world use cases due to limitations
in CPU hours and/or memory. At this point, curve fitting comes in handy. Following
Ref.~\cite{Bertsekas:2019a}, we approximate the optimal value function in the Bellman
equation for each time step separately by a linear combination of standardized
functions. This is a promising approach because the optimal value function is expected to
be continuous in the state variable. More precisely, for each $t$ we write
\begin{equation}
	\label{eq:speed_profile_12}
	\valFuncOptAppr{t}(\mystate{t};c) = \sum_{j=0}^{m-1} c_{j}^{(t)} \phi_{j}(\mystate{t})
\end{equation}
with $m \in \setN$, and $\phi_{j} \colon \setStates{t} \to \setR$ denote suitable standardized
(basis) functions. Here, we used the circumflex notation to indicate that the optimal value
function is approximated, and the superscript of the expansion coefficients emphasizes that
there is a separate approximation for each time step. For the sake of brevity, we omit the dependence
of $\valFuncOptAppr{t}$ on the coefficient vector $c$ for the remainder of this section.
The value iteration takes the same form as before (see \cref{speed_profile_general}) but with
the subtle yet important difference that the optimal value function is approximated at every
time step. Starting from the end of the horizon and moving backward in time, the approximation
at time step $t$ is then used for $t - 1$. We start at the finite horizon $t=T$ with the boundary
condition $\valFuncOptAppr{T}(\mystate{T}) = D(\mystate{T})$ and then proceed backward with
\begin{equation}
	\label{eq:speed_profile_13}
	\valFuncOptAppr{t}(\mystate{t}) = \min_{\action{t} \in \setActions{t}}
	\bigl\{ C_{t}(\mystate{t}, \action{t}) + \valFuncOptAppr{t+1}(\mystate{t+1}) \bigr\}
\end{equation}
for $t = T-1, \dotsc, 0$. Here, it is worth to empasize that at every time step, this value iteration
uses only a subset of the state space due to the involved approximation,
which allows for an efficient treatment of large state spaces. In the literature, this idea is also
known as fitted value iteration \cite{Bertsekas:2019a}. Eventually, the resulting sequence
$\valFuncOptAppr{T}, \dotsc, \valFuncOptAppr{0}$ of optimal approximated value functions is used
in the usual policy iteration \labelcref{eq:speed_profile_05}. Since these optimal value functions
are now available for (in principle) arbitrary state space variables thanks to the continuous
(yet fitted) nature of \cref{eq:speed_profile_12}, the restrictions for the discretization
of the action space are greatly relaxed.

For example, here in our specific setup, a rather fine-grained discretization of the action space
could be employed in order to get close to the analytical solution. Of course, the above explained
technique of a fitted value iteration applies generally to many dynamic programming problems
\cite{Bertsekas:2019a} and is not restricted to our problem of an optimized speed profile.

Now, we would like to present our numercial results. We solve the Bellman equation with its
ingredients given in Eqs.~(\ref{eq:speed_profile_06})--(\ref{eq:speed_profile_09}).
Again, we neglect currents and set $\ell = 100$, $\vmax = 50$, and $T = 4$. The analytical solution
according to \cref{eq:speed_profile_11} is therefore given by
$\action{0}^{\opt} = \dotso = \action{3}^{\opt} = 25 \+ \alpha / (\alpha + 1)$, which is the optimal
solution for an arbitrary $\alpha < \infty$. The Bellman equation is solved grid-based as well as with
approximating the optimal value function in the sense of a fitted value iteration. For the latter, the
approximation task---which is in this case nothing but curve fitting---is accomplished classically
as well as using the QUBO formulation of \cref{curve_fitting_qubo} solved on a quantum annealer.%
\footnote{Again, the \textit{Advantage2} QPU but with 50 runs, 256 reads per run, and an annealing time
of 1000 microseconds.} 
We employ triangular functions for the expansion \labelcref{eq:speed_profile_12},
$\phi_{j} = \triang_{j}$ (see \cref{eq:triangular}), with $m = 9$, $d = 9$, and $p = d - 1 = 8$
for the QUBO.

\begin{table*}[t]%
	\centering%
	\caption{\label{tab:solution_comparison}%
		Optimal policies and their respective total costs for different solution strategies
		with penalty parameter $\alpha = 100$, $m = 9$, and $d = 9$.
	}%
	\vspace*{0.25em}%
	\begin{ruledtabular}%
	\begin{tabular}{lccccc}%
		& $u_{0}^{\opt}$ & $u_{1}^{\opt}$ & $u_{2}^{\opt}$ & $u_{3}^{\opt}$ & $\valFuncOpt{0}(\mystate{0})$ \\%
		\cmidrule{2-6}%
		analytic & $24.752$ & $2.752$ & $24.752$ & $24.752$ & $1.990099$ \\%
		grid & $ 23.470$ & $24.490$ & $24.490$ & $26.492$ & $1.992070$ \\%
		inverse & $25 $ & $25$ & $25$ & $24.039$ & $1.990384$ \\%
		tabu & $25 $ & $25$ & $25$ & $24.039$ & $1.990384$ \\%
	   	qa & $12.5 $ & $12.5$ & $12.5$ & $50$ & $2.597069$ \\%
	\end{tabular}%
	\end{ruledtabular}%
	\vspace*{0.5em}%
\end{table*}

\begin{figure*}[t]
	\centering%
	\includegraphics[scale=0.42]{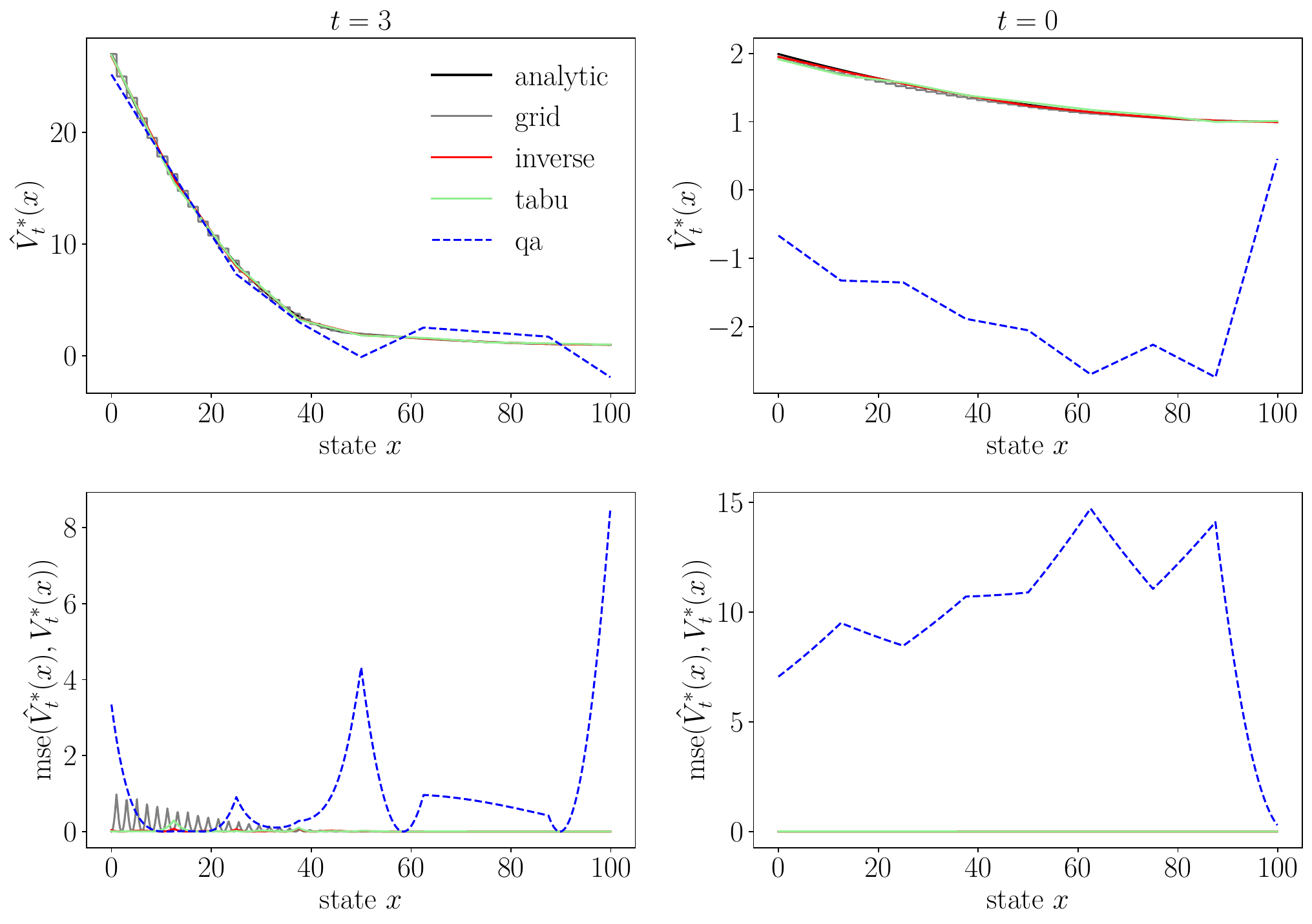}%
	\vspace{-0.75em}%
	\caption{\label{fig:value_function_n50_m9_p8}%
		Visualization of the value function approximation using different solution strategies (upper row)
		and their corresponding mean squared errors (MSE) (lower row). The left column shows the results
		for time slice $t=3$, while the right column presents the same analysis for the first time slice
		$t=0$. The upper panels depict the approximated value functions, whereas the lower panels show the
		associated approximation errors. Here and in the following figures, we use the MSE in contrast to
		the RMSE as defined in \cref{eq:sample_data_04}.
	}%
\end{figure*}

In the following bullet points, we summarize our solution strategies for the Bellman equation. The last three
methods use a fitted value iteration as described above, where the approximation task is solved classically
by matrix inversion or via a QUBO using a tabu search as well as quantum annealing. In particular, we have:
\begin{itemize}[leftmargin=*,label=\textbullet]
\item \textit{Analytic:} Approaching the problem in terms of the value function
\labelcref{eq:speed_profile_09}, we find the analytical solution from $\grad_{\? u} \valFunc{s} = 0$.
\item \textit{Grid-based:} This baseline method avoids functional approximation entirely. Instead, the state
space is discretized on a equidistant grid and estimation of the remaining costs-to-go is performed directly
on grid points. It represents a brute-force approach, useful for benchmarking, but it is computationally
expensive in high-dimensional state spaces.
\item \textit{Classically (matrix inversion):} This method relies on classical linear algebra techniques.
The system of equations resulting from the linear basis function ansatz is solved using Gaussian elimination.
It serves as a deterministic baseline and provides exact coefficients.
\item \textit{Tabu search:} The QUBO formulation of the approximation problem is solved by a tabu search,
a metaheuristic method well-suited for combinatorial landscapes. It offers a flexible alternative when exact
methods are computationally unfeasible.
\item \textit{Quantum annealing:} The QUBO formulation of the approximation problem is submitted to a D-Wave
quantum annealer. This explores the potential of currently available quantum annealing devices.
\end{itemize}

\begin{figure*}[t]
	\centering%
	\includegraphics[scale=0.38]{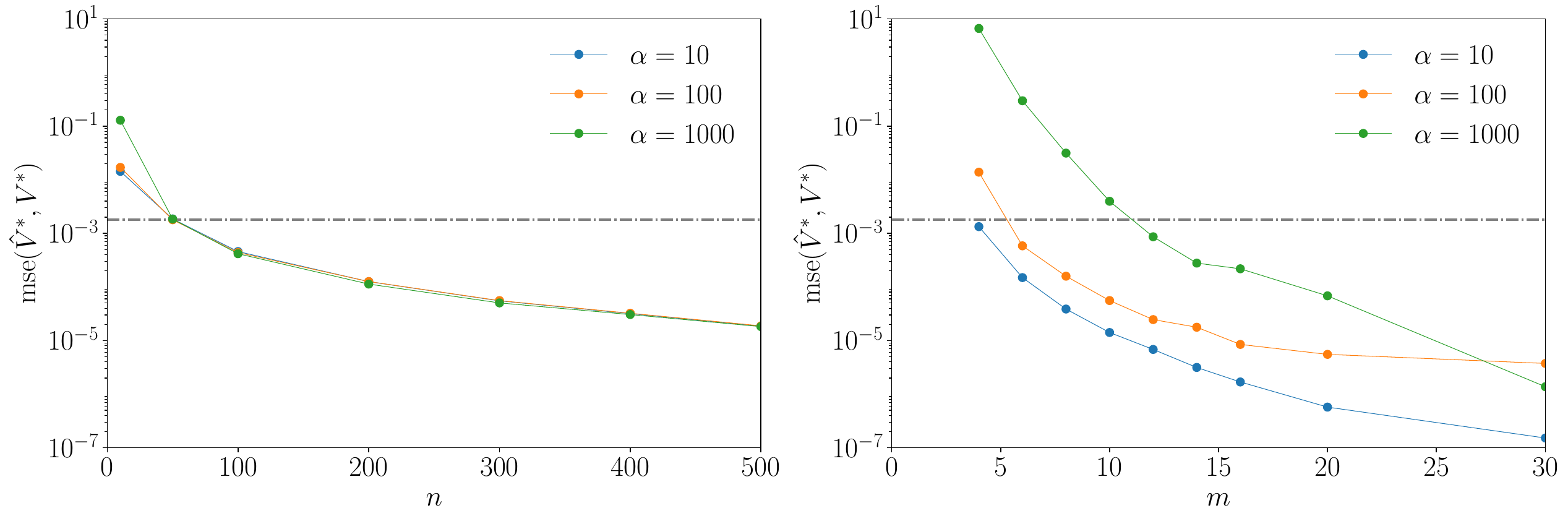}%
	\vspace*{-0.5em}%
	\caption{\label{fig:accuracy_value-function-error}%
		MSE of the reconstructed solution $\hat{V}^\opt$ compared to the reference solution $V^\opt$.
		They carry no time index because we averaged over all time slices.
		Left: Convergence behavior of the MSE with respect to the number of discretization points $n$
		of the state space for different regularization parameters $\alpha \in \{10, 100, 1000\}$.
		Right: Error decay with respect to the number of basis functions $m$ for fixed discretization $n = 50$
		and varying $\alpha$. The horizontal dash-dotted line indicates the target error threshold, where the
		approximation solution with respect to the number of basis functions $m$ starts to be more accurate
		than the grid solution.
	}%
	\vspace*{1.5em}%
	\includegraphics[scale=0.38]{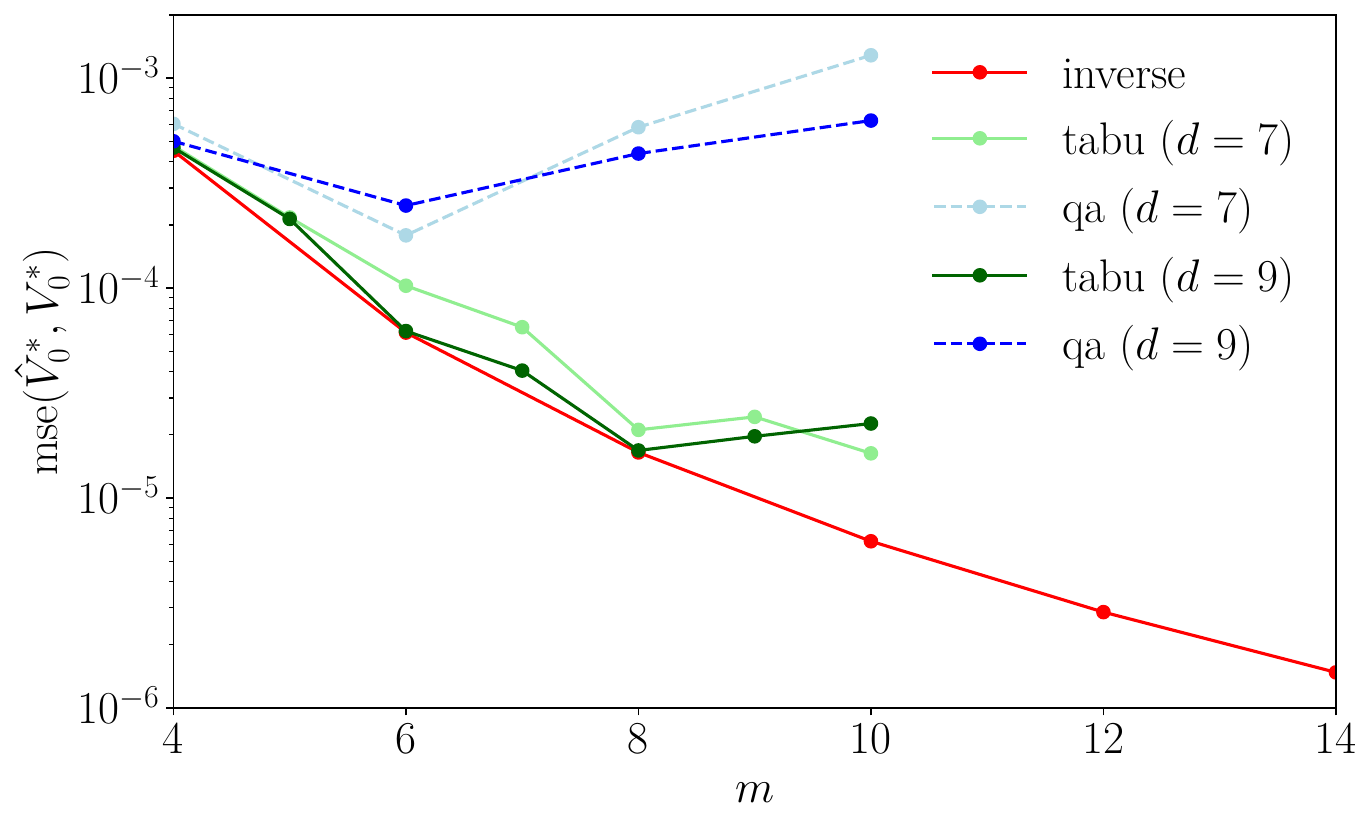}%
    \includegraphics[scale=0.38]{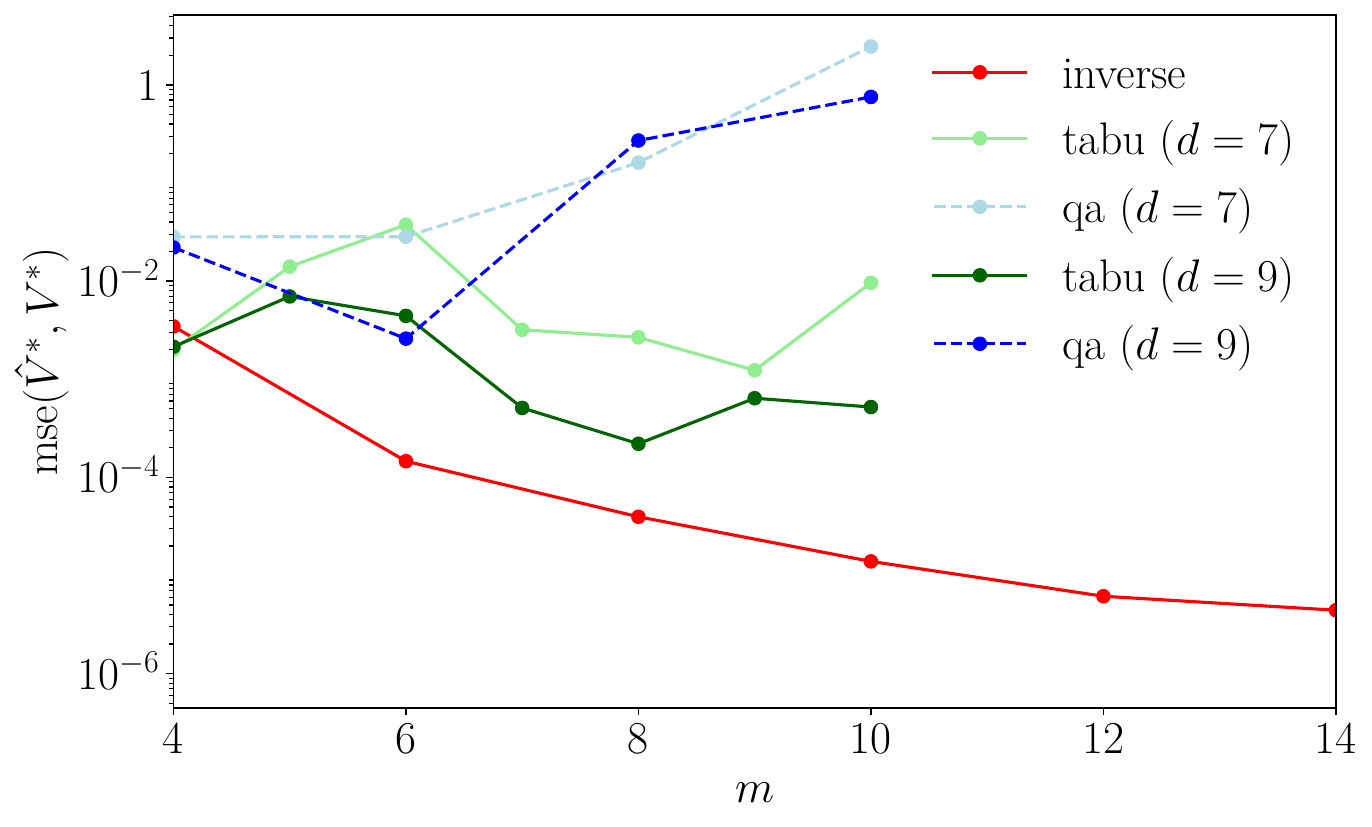}%
    \vspace*{-0.5em}%
	\caption{\label{fig:m_tabu_vs_qa}%
		Error comparison of the two solution approaches, quantum annealing and tabu search, at different
		precision levels ($d = 7$ and $d = 9$), as a function of the number of basis functions. The left
		panel shows the results for the first time slice in the fitted value iteration process, while the
		right panel displays the average performance over all time slices.
	}%
\end{figure*}

In \cref{tab:solution_comparison}, we display our results for the optimal policies obtain by means
of the different solution strategies. We have an example where the grid-based solution is indeed less optimal
than the results that employ the fitted value iteration with the approximation solved classically by
matrix inversion and QUBO-based via a tabu search. This is caused, as explained earlier, by the
discretization and forcing off-grid values onto on-grid ones. Moreover, using a QUBO solved on a quantum
annealer for the curve fitting in course of the fitted value iteration yields a worse optimal policy.

It should be noted that the choice of $m$ and $d$ results in an $81$-dimensional QUBO problem.
The annealer does not solve this problem size accurately enough, which is consistent with our findings in
\cref{numerical_experiments}. Figure \ref{fig:value_function_n50_m9_p8} shows that the solution of the
quantum annealing already exhibits a significant error at the very first time slice, preventing the method
from reaching the optimal solution. Due to error propagation, this issue is further exacerbated in the
subsequent time slices. The reasons lie in the size of the search space and the specific energy landscape
formed by the binary representation of the coefficients.

The example calculated above does not provide a clear picture of the behavior if one varies the number of
triangular functions $m$. Therefore, we varied the latter and analyzed the resulting errors. As shown in
\cref{fig:accuracy_value-function-error}, increasing the number of basis functions leads to improved
approximation results. This holds also for other values of $\alpha$. At this point, all calculations shown
in \cref{fig:accuracy_value-function-error} are classical solutions of the QUBO formulation of the
approximation task for $\valFuncOptAppr{t}$ using Gaussian elimination. Basically, this procedure corresponds
to computing a matrix inverse. The resulting solution can be interpreted as the best possible outcome within
the span of the chosen basis. In a following figure, we further investigate how different solvers approximate
this best solution. In summary, increasing the number of basis functions effectively reduces the approximation
error. The remaining challenge is to ensure that the solvers are capable of achieving this improved level of
solution quality.

Our analysis shows in \cref{fig:m_tabu_vs_qa} that as the number of basis functions---and thus the
problem size---increases, the solutions generally deviate further from the best possible solution,
referred to as the inverse solution (indicated in red). This trend holds for the tabu search but is
especially pronounced in the case of quantum annealing. Starting from six basis functions, the problem
size becomes large enough for the approximation error to increase significantly. For smaller problem
sizes, the quantum annealing yields approximation results that are comparable to those obtained with
the classical solution or the tabu search. However, beyond a certain problem size, this similarity
breaks down as clearly illustrated in \cref{fig:m_tabu_vs_qa}.

\subsection{Further discussion}
\label{further_discussion}

The proposed QUBO-based approach relates to a broader class of modern methods developed to approximate
solutions of high-dimensional HJB equations, which typically arise in optimal control and stochastic
decision problems. More specifically, we investigate a problem in the context of approximate dynamic
programming with a bounded horizon. We do not consider broader reinforcement learning techniques that
incorporate policy optimization or explicit exploration-exploitation mechanisms. Furthermore, no
explicit error-control strategies are employed, as our focus lies solely on the approximation quality
of the value function using the QUBO approach. A central challenge in solving such equations numerically
is the curse of dimensionality, i.e., the exponential growth of computational complexity with respect
to the state and action space dimensions. The QUBO-based approach can be extended to multidimensional
settings by employing, for example, a tensor product of basis functions.

Recent literature has addressed this challenge using learning-based strategies, most notably through the
approximation of the value function or the HJB operator via deep neural networks. Notable examples include
the Deep Galerkin Method \cite{Sirignano:2018} and physics-informed neural networks
\cite{Raissi:2017,Raissi:2019}. Both utilize neural networks as function approximators trained on
residuals of partial differential equations. These methods are well-suited for high-dimensional spaces by
leveraging the expressive capacity and scalability of deep architectures.

In contrast to these continuous, gradient-based approaches, our method introduces a discrete, combinatorial
framework for function approximation based on a QUBO. This discrete optimization perspective can be
interpreted as a learning approach that bypasses gradient descent and instead exploits a new perspective
on value function approximation, particularly in problem settings that can be nonconvex or nondifferentiable.
In such cases, the resulting value function may exhibit sharp transitions or multiple local extrema.
When approximating such functions using triangular  functions on a discretized domain, the associated
QUBO landscape can become irregular, potentially containing many local minima. This may lead to
instabilities in the QUBO optimization and reduced approximation accuracy. However, this is not an
intrinsic limitation of the QUBO formulation itself, but rather a structural artifact of the discretization
and basis choice. In fact, quantum annealing metaheuristics, such as those employed by the D-Wave,
offer the potential to escape local minima via quantum tunneling, which may be advantageous in
overcoming the challenges of nonconvexity.

Therefore, while nonconvexity and nondifferentiability may degrade the approximation landscape, the global
nature of the QUBO formulation together with carefully chosen basis functions preserves the possibility of
capturing the correct global structure of the value function.

Finally, we would like to note that the practical implementation of the speed profile optimization as
detailed in \cref{speed_profile_example} is based on a simplified scenario.
This simplification is intentional, as the full optimization problem rapidly becomes computationally
expensive due to the curse of dimensionality. Therefore, we deliberately considered a reduced,
one-dimensional model to illustrate the dynamic programming approach. Despite its simplicity,
our model already incorporates an important factor, the magnitude of the projected sea current,
which directly affects the optimal speed profile and fuel consumption. Moreover, the proposed framework
is flexible and can be readily extended to include additional physical effects. For example, it could be
extended to account for wind and wave conditions, variations depending on ship draft or trim, or even
to compute an optimal acceleration profile instead of a speed-only control. These extensions would
allow a more realistic representation of maritime operations while maintaining the systematic
structure provided by the dynamic programming approach. Such a comprehensive study, however, is beyond
the scope of the present work and is left for future research.

\section{Summary and conclusions}
\label{summary}

Finding optimal speed profiles is an optimization problem that grows exponentially in time.
Dynamic programming tames the combinatorial explosion of trajectories by solving the problem
backward in time. The state-dependent optimal value function as a solution of the Bellman equation
solves the problem linearly in the number of time steps but it requires a grid-based exploration
of the entire state space. If the dimension of the latter is high, the well-known phenomena of the
curse of dimensionality poses a severe obstacle, particularly for large real-world use cases.
An approximation of the optimal value function and thus its continuous representation mitigates
the rounding error that occurs in the discrete grid-based setup.

Such approximations can often be reduced to curve fitting in one dimension. In this work,
we formulated least-squares curve fitting, where we allowed for a rather general approximator
in form of a finite linear combination of standardized (basis) functions, as a QUBO that is
suited to be solved with quantum annealing. For simple functions/data points, we found that
quantum annealing in its current state is in principle able to deliver comparable results
to a classical computer. However, curve fitting tasks that result in a QUBO size larger than roughly about
$60 \times 60$ (i.e., with about 60 logical variables or more) cannot be solved accurate enough on a
D-Wave annealer yet. The reason for that is the difficult low-energy landscape of the resulting Ising
Hamiltonian. Moreover, QUBO problems described by global basis functions (e.g., Chebyshev polynomials)
are even more difficult to solve, such that even smaller problems can have poor quality.

In general, optimization problems formulated in terms of a QUBO can be modified by fine-tuning
the penalty coefficients of the constraints in order to manipulate these landscapes in a proper
way \cite{Roch:2022}. In our case, however, such an approach is not possible because the values
of the QUBO matrix entries are solely fixed by the to-be-approximated data and the
basis functions used for the approximator. Furthermore, the embedding of a floating point number
on a quantum annealer might also contribute to an already complicated low-energy landscape.

Though quantum annealing devices promise a large number of qubits, attention should be paid to
how many direct connections exist between the physical qubits. If one needs many interconnected
variables, logical qubits must be created from several physical ones. We found that basis functions
with support on the whole real line, e.g., all polynomials, result in a highly-connected QUBO
matrix that is difficult to embed into the working graph of the annealer. However, we are confident
that this issue will vanish eventually in the near future because of the advances in the ongoing
development of quantum annealing hardware.

Nevertheless, the quantum annealing-oriented curve fitting principally works. Regarding a real-world
use case---albeit still in an exploratory fashion---we investigated the linear-quadratic
optimal control problem of just-in-time arrival of a vessel, where we approximated
the value function in the Bellman equation. With an increasing number of basis functions, the solution
space can be better represented, leading to more accurate value estimates compared to the grid-based
method. There was hope that the quantum annealing process would yield lower errors in larger search
spaces. Unfortunately, we find that D-Wave's current quantum annealing hardware does not yet appear to be
capable of meeting this expectation.

Finally, we note that the present work focuses on analyzing the accuracy
and solution quality of the quantum annealing approach to curve fitting and speed profile optimization
under well-defined and controlled conditions. For this reason, potential performance enhancements
through QUBO preprocessing techniques (e.g., sparsification) and comparisons with quantum-inspired
classical solvers like SimCIM or NGQ are not considered at this stage. Runtime optimization and scaling
behavior represent important directions for future research. Such investigations will help to clarify
whether the current limitations of quantum annealing are intrinsic to the method itself or primarily
related to currently existing hardware constraints. By first establishing a reliable accuracy baseline,
this study lays the groundwork for us for subsequent studies on quantum performance, efficiency, and
practical applicability.

\vfill

\begingroup
\small
\textit{Acknowledgments.}\quad
We acknowledge funding and support from the German Federal Ministry
for Economic Affairs and Climate Action (Bun\-des\-mi\-nis\-te\-ri\-um
f\"{u}r Wirt\-schaft und Kli\-ma\-schutz) through the PlanQK initiative
(01MK20005).
This project (HA project no.~1362/22-67) is financed
with funds of LOEWE---Lan\-des-Of\-fen\-si\-ve zur Ent\-wick\-lung
Wis\-sen\-schaft\-lich-\"{o}ko\-nom\-ischer Ex\-zel\-lenz,
För\-der\-linie 3:~KMU-Ver\-bund\-vor\-ha\-ben (State Offensive for
the Development of Scientific and Economic Excellence).

\textit{Author contributions.}\quad
Conceptualization, funding acquisition, and supervision: WM, BH.
Formal analysis, investigation, and software: PI, DJ.
Validation: PI, DJ, FP.
Writing (original draft): PI, DJ.
Writing (review and editing): PI, DJ, FP.

\textit{Competing interests.}\quad
The authors have no competing interests to declare that are relevant to this article.
\endgroup


\urlstyle{same}
\bibliography{CurveFittingQUBO}

\end{document}